\documentclass[12pt,a4paper, twoside]{article}

\usepackage{amsmath,amssymb,amsthm,amsfonts,mathrsfs,amscd,environ}
\usepackage{dsfont}
\usepackage{latexsym,enumerate,color,geometry,extarrows}

\usepackage{verbatim,fancyhdr,enumitem}
 \NewEnviron{ews}{%
\begin{equation}\begin{split}
  \BODY
\end{split}\end{equation}
}

\NewEnviron{ews*}{%
\begin{equation*}\begin{split}
  \BODY
\end{split}\end{equation*}
}

\newtheorem{prop}{Proposition}[section]
 \newtheorem{thm}{Theorem}[section]
 \newtheorem{lem}[thm]{Lemma}
 
 \newtheorem{exa}[thm]{Example}
 
 \newtheorem{defn}{Definition}[section]
 \newtheorem{rem}{Remark}[section]
 
 \numberwithin{equation}{section}

\def\dif{{\mathord{{\rm d}}}}
\def\no{\nonumber}

\def\mR{{\mathbb R}}
\def\mE{{\mathbb E}}

\def\mN{{\mathbb N}}

\def\mE{{\mathbb E}}
\def\mF{{\mathbb F}}

\def\mN{{\mathbb N}}

\def\mR{{\mathbb R}}

\def\mV{{\mathbb V}}

\def\sF{{\mathscr F}}
\def\sA{{\mathscr A}}
\def\sB{{\mathscr B}}

\def\sF{{\mathscr F}}

\def\sO{{\mathscr O}}
\def\sP{{\mathscr P}}

\def\sU{{\mathscr U}}

\def\bd{\begin{defn}}
\def\ed{\end{defn}}
\def\bp{\begin{prop}}
\def\ep{\end{prop}}
\def\bl{\begin{lem}}
\def\el{\end{lem}}
\def\bt{\begin{thm}}
\def\et{\end{thm}}
\def\br{\begin{rem}}
\def\er{\end{rem}}

\allowdisplaybreaks

\title{{\bf  McKean-Vlasov  multivalued stochastic differential equations with oblique subgradients   and related stochastic control problems }
}
\author{
{\bf Hao Wu$^{a)}$, Junhao Hu$^{a)}$,  Chenggui Yuan$^{b)}$}\\
\footnotesize{$^{a)}$School of Mathematics and Statistics,
South-Central University For Nationalities}\\
\footnotesize{ Wuhan, Hubei 430000, P.R.China}\\
\footnotesize{Email: wuhaomoonsky@163.com},
\footnotesize{ junhaohu74@163.com}\\
\footnotesize{$^{b)}$Department of Mathematics, Swansea University, Bay Campus, Swansea, SA1 8EN, UK}\\
\footnotesize{Email: C.Yuan@Swansea.ac.uk}
}

\begin{document}

\maketitle

\begin{abstract}
In this article, we prove the existence   of weak solutions as well as  the existence and  uniqueness of strong solutions for McKean-Vlasov multivalued stochastic differential equations with oblique subgradients   (MVMSDEswOS, for short) by means of the  equations of Euler type and Skorohod's representation theorem.  For this type of equation, compared with the method in \cite{RWH, GRR}, since we can't  use the maximal monotony property of its constituent subdifferential operator,  some different specific techniques are applied to solve our problems.  Afterwards, we give an example for  MVMSDEswOS with time-dependent convex constraints, which can be reduced to MVMSDEswOS. Finally, we consider an optimal control problem and establish the dynamic programming principle for the value function.
\end{abstract}\noindent

\noindent


\section{Introduction}
As is well known, multivalued stochastic differential equations(MSDEs)    are widely applied to model stochastic systems  in  different branches of science and industry.  Stability, boundedness and applications    of the solution   are  the most pop research topic in the field of stochastic dynamic systems and control.     The form of these equations are as follows:
\begin{equation}\label{a}
\begin{cases}
  &X(t)+ A(X(t))\dif t\ni f( X(t))\dif t +g(   X(t))\dif  W(t)\dif t,  t>t_{0}\\
 &X(t_{0})=x_{0}\in \overline{D(A)},
 \end{cases}
 \end{equation}
where $A$ is a multivalued maximal monotone operator on $\mR^{m},$  $\overline{D(A)}$ is the closure of $D(A)$ and  $W$ is a standard Brownian motion. Compared with the usual stochastic differential equations (SDEs), i.e., $A=0,$  most of difficulties for MSDEs come from the presence of  the finite-variation process $\{k(t), t\in [0,T]\}.$ One only knows that $\{k(t), t\in [0,T]\}$ is continuous process with finite total variation and can not prove any further regularity.  This type of equations was first studied by C\'{e}pa \cite{C1,C2}. Later,
Zhang \cite{Z1} extended C\'{e}pa's results to the infinite-dimensional case.  Ren et al. \cite{RXZ} studied the large deviations for MSDEs which
solved moderate deviation problem  for the above equations.   In particular, if $A$  is taken as some subdifferential operator, the corresponding MSDEs can be used to solve a class of stochastic differential equations with reflecting boundary conditions. For MSDEs, more applications can be found in \cite{RW, RWH, PR,EE} and  references therein.

 Recently, Gassous et al. \cite{GRR} built a fundamental framework for the following  MSDEs with generalized subgradients:
\begin{align*}
\begin{cases}
 &\dif x(t)+ H(x(t))\partial \Pi(x(t))\dif t\ni f(x(t))\dif t + g(  x(t))\dif  B(t), t\in [t_{0}, T],\\
 & x(t)=x_{0}, t= t_{0},
 \end{cases}
 \end{align*}
where the new
quantity $H(\cdot, \cdot): \Omega\times \mR^{m}\rightarrow \mR^{m\times m}$  acts on the set of subgradients and $\partial\Pi(\cdot)$ is a subdifferential operator. The product $H(x(t)) \partial \Pi(x(t))$ will
be called, from now on, the set of oblique subgradients. The problem becomes challenging
due to the presence of this new term, since this new term preserves neither the monotonicity of the subdifferential operator
nor the Lipschitz property of the matrix involved.   Some different specific techniques have been applied to solve the existence and uniqueness results for this type of equations.  Later,  Gassous et al. \cite{GRR1} and Maticiuc and  Rotenstein \cite{MR+} investigated backward stochastic differential equations (BSDEs)  with oblique subgradients.   But there are few applications of MSDEs with  oblique subgradients.

  On the other hand, many researchers are interested in  studying  the following equations, which are called Mckean-Vlasov stochastic differential equations (MVSDEs):
\begin{equation*}
\begin{cases}
  &x(t)= x_{0}+\int^{t}_{0}f( x(s), \mu_{s})\dif s +g(   x(s), \mu_{s})\dif  B(s), t\in [t_{0}, \infty),\\
 &\mu_{t}:= \mbox{the probability distribution of } x(t).
 \end{cases}
 \end{equation*}
Obviously, the coefficients involved depend not only on  the state process but also on its distribution.  MVSDEs, being clearly more involved than It\^{o}'s SDEs, arise in
 McKean \cite {M2}, who was inspired by Kac's Programme in Kinetic Theory \cite{K}, as well as in some other areas of high interest such as
propagation of chaos phenomenon, PDEs, stability, invariant probability measures, social science, economics, engineering, etc. (see  e.g. \cite{CD, MH, MMLM, Z1, GA,S, W, GA, DQ1, DQ2, P}).

Motivated by the above articles, we shall study the following Mckean-Vlasov multivalued stochastic differential equations with  oblique subgradients  (MVMSDEswOS):
\begin{equation}\label{1}
\begin{cases}
 &\dif x(t)+ H(x(t), \mu_{t})\partial \Pi(x(t))\dif t\ni f(x(t), \mu_{t})\dif t  +g( x(t), \mu_{t})\dif  B(t), t\in [t_{0}, T],\\
 & x(t_{0})=x_{0},
 \end{cases}
 \end{equation}
where $H(\cdot, \cdot, \cdot): \Omega\times \mR^{m} \times \sP_{2}(\mR^{m}) \rightarrow \mR^{m\times m},$ $\mu_{t}$ is the distribution of $x(t)$ and $x_{0}\in \mR^{m}.$     The appearance of $H(\cdot, \cdot)$ leads to $H\partial \Pi$ does not inherit the maximal monotonicity of the subdifferential operator.  As a consequence, some specific techniques
when approaching the above problem are mandatory.

The main contributions of the paper are as follows:
\begin{itemize}
\item[$\bullet$]Since the coefficients involved depend not only on  the state process but also on its distribution, the method in \cite{GRR} is ineffective. We use   the  equations of Euler type and Skorohod's representaton theorem to prove the existence and uniqueness result.
\item[$\bullet$]  We present an example to illustrate our theory. The example indicates  that a class of SDEs with time-dependent constraints are equivalent to some  MVMSDEswOS.
\item[$\bullet$] In \cite{GRR}, Gassous et al. studied MSDEswOS, but they  did not give the applications for    this type of   equations.  As far as we known, there are  no works  on  optimal control problem for  MSDEswOS.  We will investigate an optimal control problem MVMSDEswOS     and establish the dynamic programming principle for the value function.
\end{itemize}

We close this part by giving our organization  in this article. In Section 2, we introduce some necessary notations,   subdifferential operators. In Section 3,  We give our main results and  an example to illustrate our theory.  In Section 4, we consider an optimal control problem and establish the dynamic programming principle for the value function.
 Furthermore,  we make the following convention: the letter $C(\eta)$ with or without indices will denote different positive constants which only depends on $\eta$, whose value may
vary from one place to another.

\section{Notations, Subdifferential operators  }
  \subsection{Notations}Throughout this paper,  let  $(\Omega, \sF, \mF, P)$  be  a complete probability space with filtration $\mF := \{\sF_{t}\}_{t\geq 0}$ satisfying the usual conditions(i.e., it is increasing and right continuous,  $\sF_{0}$ contains all $P$-null sets) taking along
 a standard  $d$-Brownian motion $B(t).$  For $x, y \in \mR^{m},$ we use $|x |$  to denote the Euclidean norm of
$x,$ and use  $\langle x, y\rangle$ to denote the Euclidean inner product. For $M\in \mR^{ m\times d},$     $|M |$ represents  $\sqrt{\mathrm{Tr} (MM^{\ast})}.$  Denote  $\sB(\mR^{d})$  by  the  Borel $\sigma-$algebra on $\mR^{d}.$  Let $\sP(\mR^{m})$ be the space  of all probability measures, and  $\sP_{p}(\mR^{m})$ denotes the space  of all probability measures defined on $\sB(\mR^{m}) $ with finite $p$th moment:
$$W_{p}(\mu):=\bigg(\bigg.\int_{\mR^{m}}|x|^{p}\mu(\dif x)\bigg)\bigg.^{\frac{1}{p}}<\infty.$$
For $ \mu, \nu\in \sP_{p}(\mR^{m})$, we define the
Wasserstein distance for $p\geq 1$  as follows: $$W_{p}(\mu, \nu):=\inf_{\pi \in \Pi(\mu, \nu)}\bigg\{\bigg. \int_{\mR^{m}\times \mR^{m}}|x-y|^{p}\pi(\dif x, \dif y)   \bigg\}\bigg.^{\frac{1}{p}}, $$
where $\Pi(\mu, \nu)$ denotes the family of all couplings for $\mu, \nu.$    Next, we  define several spaces for future use.\\

 $C([t_{0},T]; \mR^{m})$  stands for the space of all continuous functions from $[t_{0},T]$ to $\mR^{m},$ which is endowed with the uniform norm $| \varphi|_{\mathcal{C}(t_{0}, T)} = \sup_{t_{0}\leq t \leq T}|\varphi(t)|$.\\

$\mV_{0}$ denotes  the set of all continuous functions $k: [0, T]\rightarrow \mR^{m}$ with finite variation and $k(0)=0.$  $\updownarrow k\updownarrow^{t}_{s}$ stands for  the variation of $k$ on $[s,t],$ and denote $\updownarrow k\updownarrow_{t}=\updownarrow k\updownarrow^{t}_{0}.$\\

$L^{2}(0,T; \mR^{m}):=\bigg\{\bigg.  \varphi \,\mbox{is square integrable stochastic process i.e.}\, |\varphi|_{M^{2}}:=\bigg(\bigg. \mE\int^{T}_{0}|\varphi(s)|^{2}\dif s \bigg)\bigg.^{\frac{1}{2}}$
\begin{flushright}
$<\infty\bigg\}\bigg..$
\end{flushright}

$B(x,R)$  represents the  ball centered at $x$  with the radius $R$ and $\bar{B}(x,R)$  represents the closed ball centered at $x$ and with the radius $R.$

\subsection{Subdifferential operators}
We recall the definition of the subdifferential operators of a  proper lower semicontinuous convex functions $\Pi$ (l.s.c. for short) and the Moreau-Yosida approximation of the function $\Pi.$
\bd
Assume $\Pi: \mR^{m}\rightarrow (-\infty, +\infty)$  is  a  proper lower semicontinuous convex functions such that $\Pi(x)\geq \Pi(0)=0, \forall x\in \mR^{m}. $
Denote $D(\Pi)=\{x\in \mR^{m}: \Pi(x)<\infty \}.$
The set
$$\partial\Pi(x) =\{u\in \mR^{m}: \langle u, v-x   \rangle + \Pi(x)\leq \Pi(v), \forall v \in \mR^{m}  \}    $$
is called the subdifferential operator of  $\Pi.$
Denote its domain by $D(\partial\Pi)=\{x\in \mR^{m}: \partial\Pi(x)\neq \emptyset    \},$ and denote $Gr(\partial\Pi)=\{(a,b)\in \mR^{2m}:a\in \mR^{m}, b\in \partial\Pi(a)\},$  $\sA:=\bigg\{\bigg. (x,k): x \in C([0, T], \overline{D(\partial\Pi)}),  k\in \mV_{0}, \dif k(t)\in \partial \Pi(x(t)) \dif t, \mbox{and}\,\,\langle x(t)-a, \dif k(t)- b\dif t \rangle \geq 0 \,\,\mbox{for any}\,\,(a,b)\in Gr(\partial \Pi)  \bigg\}\bigg..$

\ed

The corresponding Moreau-Yosida approximation of $\Pi$    is defined as follows:
$$\Pi_{\epsilon}(x)=\inf\left\{\frac{1}{2\epsilon}|z-x|^{2}+\Pi(z): z\in \mR^{d}\right\}.$$
We know that  $\Pi_{\epsilon}$ is a convex $C^{1}-$class function.
For any  $x\in \mR^{m},$ denote $J_{\epsilon}x=x-\epsilon\nabla \Pi_{\epsilon}(x),$  where $\nabla$ is gradient operator.     Then we have  $\Pi_{\epsilon}(x)=\frac{1}{2\epsilon}|x-J_{\epsilon}x|^{2} +\Pi(J_{\epsilon}x).$
  We present  some  useful properties on the above approximation tools (for more details, see, e.g., \cite{MR+}).
\begin{itemize}
\item[(a)] $\Pi_{\epsilon}(x)=\frac{\epsilon}{2}|\nabla \Pi_{\epsilon}(x)|^{2}+\Pi(J_{\epsilon}x).$
\item[(b)] $\nabla \Pi_{\epsilon}(x)\in \partial \Pi(J_{\epsilon}x).$
\item[(c)] $|\nabla \Pi_{\epsilon}(x)- \nabla \Pi_{\epsilon}(y) |\leq \frac{1}{\epsilon}|x-y|.$
\item[(d)] $\langle \nabla\Pi_{\epsilon}(x)- \nabla \Pi_{\epsilon}(y), x-y  \rangle \geq 0.$
\item[(e)] $\langle \nabla\Pi_{\epsilon}(x)- \nabla \Pi_{\epsilon'}(y), x-y  \rangle \geq -(\epsilon +\epsilon')\langle  \nabla\Pi_{\epsilon}(x), \nabla \Pi_{\epsilon'}(y)      \rangle.$
\item[(f)] $0=\Pi_{\epsilon}(0)\leq \Pi_{\epsilon}(x), J_{\epsilon}(0)= \nabla\Pi_{\epsilon}(0)=0.$

\item [(g)]$\frac{\epsilon}{2}|\nabla\Pi_{\epsilon}(x)|^{2}\leq \Pi_{\epsilon}(x)\leq \langle \nabla\Pi_{\epsilon}(x), x  \rangle, \forall x\in \mR^{d}.$
\end{itemize}
\bl\label{ab}
Suppose that
 $Int(D(\partial\Pi))\neq \emptyset,$ where $Int(D(\partial \Pi))$ denotes the interior of  $D(\partial \Pi).$  Then for any $a \in Int(D(\partial \Pi)), $  there exists constants $\lambda_{1}> 0,\lambda_{2}\geq 0,\lambda_{3}\geq 0$ such that for any $(x, k)\in \sA $ and $0\leq s\leq t\leq T,$
$$\int^{t}_{s}\langle x(r)- a,   \dif k(r)\rangle \geq  \lambda_{1}\updownarrow k\updownarrow^{t}_{s} - \lambda_{2}\int^{t}_{s}|x(r)-a|\dif r- \lambda_{3}(t-s). $$
\el

\subsection{Strong and weak solutions of Eq.\eqref{1}}
\bd
A pair of continuous  processes $(x,k)$ is called a strong solution of \eqref{1} if
\begin{itemize}
\item[\rm{(i)}] $P(x(t_0)=x_{0})=1;$
\item[\rm{(ii)}] $x(t)$  is $\sF_{t}-$adapted.
\item[\rm{(iii)}] $(x, k)\in \sA, a.s., P.$

\item[\rm{(iv)}] $\int^{T}_{t_0}(|f(x(t), \mu_{t})|+|g(x(t), \mu_{t})|^{2})\dif t< +\infty, a.s.,P.$
\item[\rm{(v)}] For $y\in C([t_0, +\infty); \mR^{m})$ and $t_{0}\leq s \leq t < \infty$, it holds that
\begin{align*}
\int_{s}^{t}\langle y(r)-x(r),   \dif k(r)\rangle +\int_{s}^{t} \Pi(x(r))\dif r \leq \int_{s}^{t}\Pi (y(r))\dif r.
\end{align*}
\item[\rm{(vi)}] $(x,k)$ satisfies the following equation
$$x(t)+ \int^{t}_{t_{0}}H(x(s), \mu_{s})\dif  k(s)= x_{0}+\int^{t}_{t_0}f( x(s), \mu_{s})\dif s +\int^{t}_{t_0}g(   x(s), \mu_{s})\dif  B(s), t\in [t_0, T].$$
\end{itemize}
\ed
\bd
We say that Eq.\eqref{1} admits a weak solution if there exists a filtered probability space $(\tilde{\Omega}, \tilde{\sF}, \{\tilde{\sF}_{t}\}_{t\in [0,T]}, \tilde{P})$ taking along
 a  standard Brownian motion $\tilde{B}(t)$ as well as  a pair of continuous  processes $(\tilde{x},\tilde{k})$ defined on $(\tilde{\Omega}, \tilde{\sF}, \{\tilde{\sF}_{t}\}_{t\in [t_0,T]}, \tilde{P})$  such that
\begin{itemize}
\item[\rm{(i)}] $\tilde{P}(\tilde{x}(t_0)=x_{0})=1;$
\item[\rm{(ii)}] $\tilde{x}(t)$  is $\tilde{\sF}_{t}-$adapted.
\item[\rm{(iii)}] $(\tilde{x}, \tilde{k})\in \sA, a.s., \tilde{P}.$
\item[\rm{(iv)}] $\int^{T}_{0}(|f(\tilde{x}(t),\mu_{t})|+|g(\tilde{x}(t), \mu_{t})|^{2})\dif t< +\infty, a.s.,\tilde{P}.$
\item[\rm{(v)}] For  $y\in C(t_0, +\infty; \mR^{m})$ and $t_{0}\leq s \leq t < \infty$, it holds that
\begin{align*}
\int^{t}_{s}\langle y(r)-\tilde{x}(r),   \dif \tilde{k}(r)\rangle +\int^{t}_{s} \Pi(\tilde{x}(r))\dif r \leq \int^{t}_{s}\Pi (\tilde{y}(r))\dif r,
\end{align*}
\item[\rm{(vi)}] $(\tilde{x},\tilde{ k})$ satisfies the following equation
$$\tilde{x}(t)+ \int^{t}_{t_{0}}H(\tilde{x}(s), \mu_{s})\dif \tilde{ k}(s)= \tilde{x}_{0}+\int^{t}_{t_0}f( \tilde{x}(s), \mu_{s})\dif s +\int^{t}_{t_0}g(\tilde{x}(s), \mu_{s})\dif  B(s), t\in [t_0, T].$$
\end{itemize}
\ed
\br
We make the following convention: For a pair of process $(x, k)$ satisfying  $\dif k(t)\in \partial \Pi(x(t)) \dif t $,  we denote that   $U(t) $ is a process such that $\dif k(t)= U(t) \dif t.$

\er

\section{Main Results}
 Before giving our main results for Eq.\eqref{1}.  For the sake of simplicity, we assume $f(0, \delta_{0} )=g(0,\delta_{0})=0,$  where $\delta_{x}$ stands for the Dirac measure at $x.$   Now,  we make the following assumptions:
\begin{itemize}
\item[(A1)]The coefficients $f, g$  satisfy that  for some positive constant $L>0$ and all $x,y \in \mR^{d}, \mu, \nu \in \sP_{2}(\mR^{m}),$
$$|f(x,\mu)-f(y,\nu)|+|g(x,\mu)-g(y,\nu)| \leq L(|x-y|+W_{2}(\mu, \nu)), $$
Furthermore, from the above assumptions,  one has
$$|f(x,\mu)|\leq L(|x|+W_{2}(\mu, \delta_{0})), |g(x,\mu)|\leq L(|x|+W_{2}(\mu, \delta_{0})). $$

\item[(A2)] $H(\cdot, \cdot)=(a_{ij}(\cdot, \cdot))_{m\times m}: \mR^{m}\times \sP_{2}(\mR^{m})\rightarrow \mR^{m\times m}$ is a continuous mapping and for any $(x, \mu)\in \mR^{m}\times \sP_{2}(\mR^{m}),$ $H(x, \mu)$ is a invertible symmetric matrix.  Moreover,   there exist two positive $a_{H}, b_{H}$ such that
    \begin{itemize}
\item[(i)] $a_{H}|u|^{2}\leq   \langle H(x, \mu)u ,u\rangle\leq b_{H}|u|^{2}, \forall u\in \mR^{m}.$
\item[(ii)]For all $x,y \in \mR^{m}, \mu, \nu \in \sP_{2}(\mR^{m}),$
$$|H(x,\mu)-H(y,\nu)| +|H^{-1}(x,\mu)-H^{-1}(y,\nu)|\leq L(|x-y|+W_{2}(\mu, \nu)), $$
where $|H(x,\mu)|:=\bigg(\bigg.\sum^{m}_{i, j=1} |a_{ij}(x,\mu)|^{2}  \bigg)\bigg.^{\frac{1}{2}}$ and  $H^{-1}(x,\mu)$ the inverse of the matrix $H(x,\mu).$
\end{itemize}
\end{itemize}

\bl \label{L2L}
Let $I$ be  an arbitrary set of indexes. For each $i\in I,$  suppose that $(\Omega^{i}, \sF^{i}, P^{i},$ $\{\sF_{t}^{i}\}_{t\geq 0}, B^{i}, x^{i}, k^{i} )$ is a weak solution of the equation
\begin{equation}\label{i+++}
\begin{cases}
 & \dif x^{i}(t)+ H_{i}\partial \Pi(x^{i}(t))\dif t\\
 &\quad \quad \quad \ni f^{i}(x^{i}(t), \mu^{i}_{t})\dif t  + g^{i}( x^{i}(t), \mu^{i}_{t})\dif  B^{i}(t), t\in [t_{0}, T],\\
 & x(t_{0})=x^{i}_{0}\in \mR^{m}.
 \end{cases}
 \end{equation}
where $f^{i}, g^{i}, H_{i}$ satisfy $(\mathrm{A}1)-(\mathrm{A}2)$ and $H^{i}$ is independent on $ x, \mu.$  If $\sup_{i}\mE[\sup_{0\leq t \leq T}|x^{i}(t)|^{2}]< \infty,$ then $(x^{i}, k^{i})_{i\in I}$ is tight in $C([t_{0}, T]; \mR^{m})\times C([t_{0}, T]; \mR^{m}).$
\el
\begin{proof}
Set $\hat{x}^{i}(t)= H_{i}^{-\frac{1}{2}}x^{i}(t).$  From \eqref{i+++},  we know that $\hat{x}^{i}(t)$ satisfies the following equation:
\begin{equation}\label{i++}
\begin{cases}
 & \dif \hat{x}^{i}(t)+ H_{i}^{\frac{1}{2}}\partial \Pi(x^{i}(t))\dif t\\
 &\quad \quad \quad \ni H_{i}^{-\frac{1}{2}}f^{i}(x^{i}(t), \mu^{i}_{t})\dif t  + H_{i}^{-\frac{1}{2}}g^{i}( x^{i}(t), \mu^{i}_{t})\dif  B^{i}(t), t\in [t_{0}, T],\\
 & x(t_{0})=x^{i}_{0}.
 \end{cases}
 \end{equation}
 When we apply  It\^{o}'s formula to $|\hat{x}^{i}(t)|^{2}$ for the above equation, $H_{i}^{\frac{1}{2}}$  will disappear. Then, we can use the maximal monotony property  of the   subdifferential operator $\partial \Pi$.  Using  the similar method in   \cite [Theorem 3.2]{Z}, we can get the desired result.
\end{proof}

\bt\label{3.2}
Let $f,g$  satisfy the following linear growth conditions:
$$|f(x,\mu)|\leq L(|x|+W_{2}(\mu, \delta_{0})), |g(x,\mu)|\leq L(|x|+W_{2}(\mu, \delta_{0})),$$  and assume $(\mathrm{A2})$ holds.
Then Eq.\eqref{1} has a weak solution.
\et
\begin{proof}
For fixed $n,$  take $t_{n}= 2^{-n}\lfloor 2^{n}t\rfloor,$  where $\lfloor z\rfloor$ denotes the integer part of a real number $z.$ In addition, denote $\mu_{t_{0}}=\delta_{x_{0}}.$  Consider the following approximation equation for $n\geq 2$:
\begin{equation}\label{h}
\begin{cases}
 &\dif x^{n}(t)+ H(x^{n-1}(t_{n}), \mu^{n-1}_{t_{n}})\partial \Pi(x^{n}(t))\dif t\\
 &\quad \quad \quad \ni f(x^{n-1}(t_{n}), \mu^{n-1}_{t_{n}})\dif t  +g( x^{n-1}(t_{n}), \mu^{n-1}_{t_{n}})\dif  B(t), t\in [t_{0}, T],\\
 & x(t_{0})=x_{0}.
 \end{cases}
 \end{equation}
By solving a deterministic Skorohod problem (see \cite{C2} for more details), we can obtain the solution of this equation  step by step. Thus, there exists $(x^{n}, k^{n})$ such that
\begin{equation}\label{i}
\begin{cases}
 & x^{n}(t)+ \int^{t}_{t_{0}}H(x^{n-1}(s_{n}), \mu^{n-1}_{s_{n}})\dif  k^{n}(s)\\
 &\quad \quad \quad = \int^{t}_{t_{0}}f(x^{n-1}(s_{n}), \mu^{n-1}_{s_{n}})\dif s  +\int^{t}_{t_{0}} g( x^{n-1}(s_{n}), \mu^{n-1}_{s_{n}})\dif  B(s), t\in [t_{0}, T],\\
 & x(t_{0})=x_{0}.
 \end{cases}
 \end{equation}
 For given $a \in Int (D(\partial \Pi)),$  set  $\hat{x}^{n}(t)=x^{n}(t)-a,    \tilde{x}^{n}(t) =H^{-\frac{1}{2}}(x^{n-1}(t_{n}), \mu^{n-1}_{t_{n}})\hat{x}^{n}(t).$
Then, $\tilde{x}^{n}$  satisfies the following equation:
\begin{align}\label{a10+}
\tilde{x}^{n}(t)&= \int^{t}_{t_{0}}H^{-\frac{1}{2}}(x^{n-1}(s_{n}), \mu^{n-1}_{s_{n}})\dif\hat{x}^{n}(s)+|\tilde{x}^{n}(t_{0})|^{2}\no\\
& =\int^{t}_{t_{0}}H^{-\frac{1}{2}}(x^{n-1}(s_{n}), \mu^{n-1}_{s_{n}})f( x^{n-1}(s_{n}), \mu^{n-1}_{s_{n}})\dif s - \int^{t}_{t_{0}}H^{\frac{1}{2}}(x^{n-1}(s_{n}), \mu^{n-1}_{s_{n}})\dif k^{n}(s)\no\\
&+ \int^{t}_{t_{0}}H^{-\frac{1}{2}}(x^{n-1}(s_{n}), \mu^{n-1}_{s_{n}})g( x^{n-1}(s_{n}), \mu^{n-1}_{s_{n}})\dif B(s)+|\tilde{x}^{n}(t_{0})|^{2}.
\end{align}
Using It\^{o}'s formula,  we have
\begin{align}\label{j}
&|\tilde{x}^{n}(t)|^{2}= |\tilde{x}^{n}(t_{0})|^{2} + 2\int^{t}_{t_0}\langle \tilde{x}^{n}(s), H^{-\frac{1}{2}}(x^{n-1}(s_{n}), \mu^{n-1}_{s_{n}})f( x^{n-1}(s_{n}), \mu^{n-1}_{s_{n}}) \rangle\dif s\no\\
& +\int^{t}_{t_0}|H^{-\frac{1}{2}}(x^{n-1}(s_{n}), \mu^{n-1}_{s_{n}})g( x^{n-1}(s_{n}), \mu^{n-1}_{s_{n}})|^{2}\dif s\no\\
& + 2\int^{t}_{t_0}\langle \tilde{x}^{n}(s), H^{-\frac{1}{2}}(x^{n-1}(s_{n}), \mu^{n-1}_{s_{n}})g( x^{n-1}(s_{n}), \mu^{n-1}_{s_{n}})\dif B(s)\rangle  \no\\
&- 2\int^{t}_{0}\langle x^{n}(s)-a,   \dif k^{n}(s)\rangle\no\\
&=: |\tilde{x}^{n}(t_0)|^{2} + \sum^{4}_{i=1}I_{i}.
\end{align}
Firstly,  we estimate $I_{1}.$
\begin{align*}
I_{1}=\int^{t}_{t_0}&\langle \tilde{x}^{n}(s), H^{-\frac{1}{2}}(x^{n-1}(s_{n}), \mu^{n-1}_{s_{n}})f( x^{n-1}(s_{n}), \mu^{n-1}_{s_{n}}) \rangle\dif s\\
&\leq 2\int^{t}_{t_0}|\tilde{x}^{n}(s)|^{2}\dif s+  a_{H}^{-1}\int^{t}_{t_0} | f( x^{n-1}(s_{n}), \mu^{n-1}_{s_{n}})|^{2}\dif s\\
& \leq (1+8a_{H}^{-1}b_{H}L^{2})\int^{t}_{t_0}\sup_{1\leq i\leq n}\sup_{t_{0}\leq r\leq s}|\tilde{x}^{i}(s)|^{2}\dif s+ 8L^{2}a^{2}a_{H}^{-1}T.
\end{align*}
Similarly,  for $I_{2},$ we have
\begin{align*}
I_{2}=\int^{t}_{t_0}&|H^{-\frac{1}{2}}(x^{n-1}(s_{n}), \mu^{n-1}_{s_{n}})g( x^{n-1}(s_{n}), \mu^{n-1}_{s_{n}})|^{2}\dif s\\
 &\leq 8a_{H}^{-1}b_{H}L^{2}\int^{t}_{t_0}\sup_{1\leq i\leq n}\sup_{t_{0}\leq r\leq s}|\tilde{x}^{i}(s)|^{2}\dif s+ 8L^{2}a^{2}a_{H}^{-1}T.
\end{align*}
 Moreover,  from the Burkholder-Davis-Gundy (BDG) inequality,  for any $l>0,$  we have
\begin{align*}
I_{3}=\mE&\sup_{t_0\leq r\leq t}\bigg|\bigg.\int^{r}_{t_{0}}\langle \tilde{x}^{n}(s), H^{-\frac{1}{2}}(x^{n-1}(s_{n}), \mu^{n-1}_{s_{n}})g( x^{n-1}(s_{n}), \mu^{n-1}_{s_{n}})\dif B(s)\rangle\bigg|\bigg.\\
& \leq 32\mE\bigg(\bigg.\int^{t}_{t_{0}}|\tilde{x}^{n}(s)|^{2} | H^{-\frac{1}{2}}(x^{n-1}(s_{n}), \mu^{n-1}_{s_{n}})g( x^{n-1}(s_{n}), \mu^{n-1}_{s_{n}})|^{2} \dif s \bigg)\bigg.^{\frac{1}{2}}\\
& \leq 32\mE\bigg[\bigg. \sup_{t_0\leq s\leq t}|\tilde{x}^{n}(s)|\bigg(\bigg.\int^{t}_{t_{0}} |H^{-\frac{1}{2}}(x^{n-1}(s_{n}), \mu^{n-1}_{s_{n}})g( x^{n-1}(s_{n}), \mu^{n-1}_{s_{n}})|^{2} \dif s \bigg)\bigg.^{\frac{1}{2}}\bigg]\bigg.\\
& \leq \frac{1}{l}\mE[\sup_{t_0\leq s\leq t}|\tilde{x}^{n}(s)|^{2}]  +32l  \mE\int^{t}_{t_{0}}|H^{-\frac{1}{2}}(x^{n-1}(s_{n}), \mu^{n-1}_{s_{n}})g( x^{n-1}(s_{n}), \mu^{n-1}_{s_{n}})|^{2}\dif s\\
&\leq  \frac{32}{l}\mE[\sup_{1\leq i\leq n}\sup_{t_{0}\leq s\leq t}|\tilde{x}^{i}(s)|^{2}]  + 256 l L^{2}a_{H}^{-1}b_{H} \mE\int^{t}_{t_{0}}\sup_{1\leq i\leq n}\sup_{t_{0}\leq r\leq s}|\tilde{x}^{i}(r)|^{2}\dif s+ 8L^{2}a^{2}a_{H}^{-1}T.
\end{align*}
Finally, we calculate $I_{4}.$  By Lemma \ref{ab}, we have
\begin{align*}
- 2\int^{t}_{0}\langle x^{n}(s)-a,   \dif k^{n}(s)\rangle&\leq  \lambda_{2}\int^{t}_{s}|x^{n}(s)-a|\dif s+ \lambda_{3}T\\
& \leq \lambda_{2}b_{H}T\int^{t}_{s}\sup_{1\leq i\leq n}\sup_{t_{0}\leq r\leq s}|\tilde{x}^{i}(r)|^{2}\dif s+ \lambda_{3}T.
\end{align*}
Take $l= 64.$ Combining the above calculations,   we have
\begin{align*}
\mE[\sup_{1\leq i\leq n}\sup_{t_{0}\leq s\leq t}|\tilde{x}^{i}(s)|^{2}]&\leq (2^{18}a_{H}^{-1}b_{H}L^{2}+2^{5}a_{H}^{-1}b_{H}L^{2}+2\lambda_{2}b_{H}T+2 )\mE\int^{t}_{t_{0}}\sup_{1\leq i\leq n}\sup_{t_{0}\leq r\leq s}|\tilde{x}^{i}(r)|^{2}\dif s\\
&+ 8L^{2}a^{2}a_{H}^{-1}T+2\lambda_{3}T+|\tilde{x}^{n}(t_{0})|^{2}.
\end{align*}
By Gronwall's inequality,  we have
$$\sup_{ n\geq 0}\mE[\sup_{t_{0}\leq s\leq t}|\tilde{x}^{n}(s)|^{2}]\leq C(a_{H},b_{H},T, \lambda_{2},\lambda_{3} ).$$
Then,
$$\sup_{ n\geq 0}\mE[\sup_{t_{0}\leq s\leq t}|x^{n}(s)|^{2}]\leq C(a_{H},b_{H},T, \lambda_{2},\lambda_{3} ).$$
Furthermore, from Lemma \ref{ab} and \eqref{j}, we derive
$$\sup_{ n\geq 0}\mE[\sup_{t_{0}\leq s\leq t}|k^{n}(s)|^{2}]\leq C(a_{H},b_{H},T, \lambda_{2},\lambda_{3} ).$$
Using the above estimates,  by Lemma \ref{L2L}, $(x^{n}, k^{n})_{n\in \mN}$ is  tight in $C([t_{0}, T]; \mR^{2m}).$
By the Prohorov theorem there exists a subsequence which we still denote by $(x^{n}, k^{n}, \updownarrow k^{n}\updownarrow, B) $ such that,  as $n\rightarrow \infty,$
$$(x^{n}, k^{n}, \updownarrow k^{n}\updownarrow, B)\Rightarrow (x, k, \updownarrow k\updownarrow, B).$$
By Skorohod   representation theorem,  we can choose a probability space $(\check{\Omega}, \check{\sF}, \check{P})$ and some quadruples $(\check{x}^{n}, \check{k}^{n}, \check{V}^{n},   \check{B}^{n})$ and $(\check{x}, \check{k}, \check{V},   \check{B})$ defined on $(\check{\Omega}, \check{\sF}, \check{P}),$ having the same laws as $(x^{n}, k^{n}, \updownarrow k^{n}\updownarrow, B)$ and $(x, k, \updownarrow k \updownarrow, B),$ respectively, such that , in $C([t_{0}, T]; \mR^{2m+1+d}),$ as $n\rightarrow \infty,$
$$(\check{x}^{n}, \check{k}^{n}, \check{V}^{n},   \check{B}^{n})\rightarrow (\check{x}, \check{k}, \check{V},   \check{B}), a.e., n\rightarrow \infty. $$
Since $(x^{n}, k^{n}, \updownarrow k^{n} \updownarrow, B)\Rightarrow (\check{x}, \check{k}, \check{V},   \check{B}),$  then by Proposition 16, in \cite{GRR}, we have that,  for all $t_{0}\leq s<t,$
\begin{align} \label{k1}
&\check{x}(t_{0})=x_{0}, \check{k}(t_{0})=0, |\check{k}|_{v}(t)-|\check{k}|_{v}(s)\leq \check{V}(t)-\check{V}(s), 0=\check{V}(t_{0})\leq\check{V}(s)\leq \check{V}(t), \check{P}-a.e.
\end{align}
Furthermore, since $\forall t_{0}\leq s<t,$
\begin{align*}
\int^{t}_{s}\Pi(x^{n}(r))\dif r \leq \int^{t}_{s}\Pi(y(r))\dif r - \int^{t}_{s}\langle y(r)-x^{n}(r), \dif k^{n}(r) \rangle,
\end{align*}
 and  \cite[Proposition 16]{GRR}, one can see that
\begin{align}\label{k2}
\int^{t}_{s}\Pi(\check{x}(r))\dif r \leq \int^{t}_{s}\Pi(y(r))\dif r - \int^{t}_{s}\langle y(r)-\check{x}(r), \dif \check{k}(r) \rangle.
\end{align}
Combining \eqref{k1} and \eqref{k2}, it infer
$$\dif \check{k}(r)\in \partial \Pi (\check{x}(r))(\dif r).  $$
Using the Lebesgue  theorem,  we derive
\begin{align*}
\check{M}^{n}(\cdot)&=x_{0}+\int^{\cdot}_{t_{0}}f(\check{x}^{n-1}(s_{n}), \mu^{n-1}_{s_{n}})\dif s  +\int^{\cdot}_{t_{0}}g( \check{x}^{n-1}(s_{n}), \mu^{n-1}_{s_{n}})\dif  B(s)\\
& \rightarrow \check{M}(\cdot)=x_{0}+\int^{\cdot}_{t_{0}}f(\check{x}(s), \mu_{s})\dif s  +\int^{\cdot}_{t_{0}}g( \check{x}(s), \mu_{s})\dif  \check{B}(s)\,\, \mbox{in}\,\,  L^{2}(t_{0}, T; \mR^{m}).
\end{align*}
By Proposition 17 in \cite{GRR}, it follows that
$$\mathcal{L}(\check{x}^{n}, \check{k}^{n},   \check{B}^{n},\check{M}^{n} )=\mathcal{L}(x^{n}, k^{n},   B^{n},M^{n} )\,\, \mbox{in}\,\, C([t_{0}, T]; \mR^{3m+d}), $$
where $\mathcal{L}(\cdot)$ is denoted by the probability law of the random variables. Since
$$ x^{n}(t)+ \int^{t}_{t_{0}}H(x^{n-1}(s_{n}), \mu^{n-1}_{s_{n}})\dif k^{n}(s)-M^{n}(t)=0, $$
we have
$$ \check{x}^{n}(t)+ \int^{t}_{t_{0}}H(\check{x}^{n-1}(s_{n}), \mu^{n-1}_{s_{n}})\dif \check{k}^{n}(s)-\check{M}^{n}(t)=0. $$
Letting $n\rightarrow \infty, $ one has
$$ \check{x}(t)+ \int^{t}_{t_{0}}H(\check{x}(s), \mu_{s})\dif \check{k}(s)-\check{M}(t)=0. $$
This means
\begin{align*}
\check{x}(t)+ \int^{t}_{t_{0}}H(\check{x}(s), \mu_{s})\dif \check{k}(s)=x_{0}+ \int^{t}_{t_{0}}f(\check{x}(s), \mu_{s})\dif s  + \int^{t}_{t_{0}}g( \check{x}(s), \mu_{s})\dif  \check{B}(s).
\end{align*}
Thus, $((\check{\Omega}, \check{\sF}, \check{P}), \sF_{t}^{\bar{B},\bar{X}},\check{x}(t), \check{k}(t), \check{B}(t) )_{t\geq t_{0}}$ is a weak solution.  The proof is  complete.

\bt
Let $(\mathrm{A}1)-(\mathrm{A}2)$ hold and $H$  is independent of $\mu.$  The Eq.\eqref{1} has a unique strong solution.
\et
 Based on Theorem \eqref {3.2},    it suffices to prove the pathwise  uniqueness. Assume that $(x, k)$ and $(x', k')$ are two solutions of \eqref{1}.  Let $Q(t):= H^{-1}(x(t))+H^{-1}(x'(t)). $  Then,  we have
$$\dif Q^{\frac{1}{2}}(t)= \dif N(t)+ \sum^{d}_{i=1}\beta_{i}(t)\dif B_{i}(t), $$
where $N$ is an  $\mR^{m\times m}-$valued bounded variation continuous stochastic process with $N(t_{0})=0,$ and $\beta_{i}, i=1,2,\cdots, d$  is an $\mR^{m\times m}-$valued stochastic process such that $\mE\int^{T}_{t_{0}}|\beta_{i}(t)|^{2}\dif t< \infty.$
Set  $\hat{x}(t)=Q^{\frac{1}{2}}(t)(x(t)-x'(t)).$
Then, $\hat{x}(t)$  satisfies the following equation:
\begin{align*}
\dif \hat{x}(t)&=[\dif Q^{\frac{1}{2}}(t)](x(t)-x'(t))+ Q^{\frac{1}{2}}(t)\dif (x(t)-x'(t))+ \sum^{d}_{i=1}\beta_{i}(t)(g(x(t), \mu_{t})-g(x'(t), \mu'_{t}))e_{i}\dif t\\
&=: \dif F(t)+ G(t)\dif B(t),
\end{align*}
where
\begin{align*}
 \dif F(t)&=[\dif N(t)]Q^{-\frac{1}{2}}(t)\hat{x}(t)+ Q^{\frac{1}{2}}(t)(H(x(t))\dif k(t)-H(x'(t))\dif k'(t) )\\
&+Q^{\frac{1}{2}}[f(x(t), \mu_{t})\dif k(t)-f(x'(t), \mu'_{t})]\dif t+  \sum^{d}_{i=1}\beta_{i}(t)(g(x(t), \mu_{t})-g(x'(t), \mu'_{t}))e_{i}\dif t,
\end{align*}
and
\begin{align*}
G(t)&=\Gamma(t)+  Q^{\frac{1}{2}}(g(x(t), \mu_{t})-g(x'(t), \mu'_{t}))e_{i},
\end{align*}
and $\Gamma(t)$ is an $\mR^{m\times d}$ matrix with the columns $\beta_{1}(t)(x(t)-x'(t)), \cdots, \beta_{i}(t)(x(t)-x'(t)).$\\
Using the properties of matrices $H, H^{-1}$ and the following relation
\begin{align*}
Q(t)(H(x'(t))\dif k'(t)-H(x(t))\dif k(t) )
&= (H^{-1}(x(t))-H^{-1}(x'(t)))[H(x'(t))\dif k'(t)\\
&+H(x(t))\dif k(t)]+ 2(\dif k'(t)-\dif k(t)),
\end{align*}
we have
\begin{align*}
\langle &\hat{x}(t),   Q^{\frac{1}{2}}(t)(H(x(t))\dif k(t)-H(x'(t))\dif k'(t) ) \rangle\\
&=\langle x(t)-x'(t), (H^{-1}(x(t))- H^{-1}(x'(t)))(H(x(t))\dif k(t)+H(x'(t))\dif k'(t))   \rangle\\
&\quad - 2\langle x(t)-x'(t), \dif k(t)- \dif k'(t)   \rangle\\
& \leq C(L, b_{H})|x(t)-x'(t)|^{2}(\dif \updownarrow k\updownarrow_t+ \dif \updownarrow k'\updownarrow_t).
\end{align*}
Thus,
$$\langle \hat{x}(t),   \dif F(t)    \rangle + \frac{1}{2}|G(t)|^{2}\dif t \leq |\hat{x}(t)|^{2}\dif V(t),  $$
where
$$\dif V(t)= C(L, b_{H})(1+ \dif \updownarrow N\updownarrow_{t}+ \dif \updownarrow k\updownarrow_{t}+ \dif \updownarrow k'\updownarrow_{t})+ C(L, b_{H})\sum^{d}_{i=1}|\beta_{i}(t)|^{2}\dif t.     $$
By    \cite[Proposition 14]{GRR},  it yields
$$\mE\frac{e^{-2V(t)}|\hat{x}(t)|^{2}}{1+2e^{-2V(t)}|\hat{x}(t)|}\leq \mE\frac{e^{-2V(0)}|\hat{x}(0)|^{2}}{1+2e^{-2V(0)}|\hat{x}(0)|}=0.     $$
Then,
$$Q^{\frac{1}{2}}(t)(x(t)-x'(t))=0,\, \mbox{for\,all},\,t\geq t_{0}.  $$
Consequently,
$$x(t)=x'(t),\, \mbox{for\,all},\,t\geq t_{0}.$$

Then, the uniqueness holds.  The proof is complete.
\end{proof}

The following examples illustrate the theories about existence and uniqueness.
\begin{exa} {\rm
Assume that $\sO$ is a closed convex subset of $\mR^{2},$ and  that $I_{\sO}$ is the indicator function of $\sO$, i.e.,
\begin{align*}
 I_{\sO}(x)=\begin{cases}
  &+\infty, x \notin \sO,\\
 &1,  x\in \sO.\\
 \end{cases}
 \end{align*}
Then,  the subdifferential operator of $I_{\sO}$ given by

\begin{align*}
\partial I_{\sO}(x)=\begin{cases}
  &\emptyset, x \notin \sO,\\
 &\{0\},  x\in \mbox{Int} (\sO),\\
 & \Lambda_{x}, x\in \partial \sO,
 \end{cases}
 \end{align*}
where  $\Lambda_{x}$ is the exterior normal cone at $x$ and $\mbox{Int} (\sO)$ is the interior of $\sO.$   For any $(x, \mu)\in \mR^{2}\times \sP_{2}(\mR^{2}), $ set

$$
H(x, \mu)=\left (
  \begin{array}{ccc}
    \sin x +5+ \cos (W_{2}(\mu, \delta_{0})) & 0 \\
    0 & e^{\cos x}+4 + (W_{2}(\mu, \delta_{0})\wedge 1)\\
  \end{array}
\right),
$$
$$f(x, \mu)=\sqrt{|x|^{2}+5}+W_{2}(\mu, \delta_{0}), g(x, \mu)=e^{1\wedge |x|}+\sin (W_{2}(\mu, \delta_{0})).$$

Obviously, $H(\cdot, \cdot), f(\cdot, \cdot), g(\cdot, \cdot)$ satisfy $\mathrm{(A1)}$ and $\mathrm{(A2)}.$  Then,  the following equation has a unique weak solution:
\begin{align*}
\begin{cases}
 &\dif x(t)+ H(x(t), \mu_{t})\partial \Pi(x(t))\dif t= f(x(t), \mu_{t})\dif t  +g( x(t), \mu_{t})\dif  B(t), t\in [t_{0}, T],\\
 & x(t_{0})=x_{0}.
 \end{cases}
 \end{align*}
Moreover,  If we take
$$
H(x)=\left (
  \begin{array}{ccc}
    \sin x +5 & 0 \\
    0 & e^{x\wedge 1}+4 + \cos x\\
  \end{array}
\right),
$$
then  following equation has a unique strong solution:
\begin{align*}
\begin{cases}
 &\dif x(t)+ H(x(t), \mu_{t})\partial \Pi(x(t))\dif t\ni f(x(t), \mu_{t})\dif t  +g( x(t), \mu_{t})\dif  B(t), t\in [t_{0}, T],\\
 & x(t_{0})=x_{0}.
 \end{cases}
 \end{align*}

}
\end{exa}

\begin{exa} {\rm
In this example, we shall show that some MVMSDEs with time-dependent constraints can be transferred to some suitable MVMSDEswOS with oblique subgradients. More precisely, consider the following MVMSDEswOS with time-dependent constraints:
\begin{equation}\label{37}
\begin{cases}
 &\dif x(t)+\partial I_{H(t)\Xi}(x(t))\dif t\ni f(x(t),  \mu_{t}\circ H(t))\dif t  +g(  x(t),   \mu_{t}\circ H(t))\dif  B(t), t\in [t_{0}, T],\\
 & x(t)=x_{0},
 \end{cases}
 \end{equation}
where $\mu_{t}$ is the distribution of $x(t),$ $\Xi\in \mR^{d}$  is a convex set and a deterministic  time-dependent matrix  $H: [0, T] \rightarrow \mR^{m\times m}$ satisfying $(\mathrm{A2})$. Assume the coefficients $f,g$ satisfy the condition $(\mathrm{A1}).$
Furthermore, for any $\sO\in \sB(\mR^{m}), H(t)(\sO):=\{H(t)a| a\in  \sO    \},  \mu\circ H(t)(\sO):=\mu(H(t)(\sO)).$  Since $\partial I_{H(t)\Xi}(y)= \partial I_{\Xi}(H^{-1}(t)y), \forall y \in D(\partial I_{{H(t)\Xi}}),$  we can easily verify that $x$ is a solution for Eq.\eqref{37} if only if $\bar{x}=H^{-1}x$ is a solution for the following GSDEs with oblique subgradients:
\begin{equation}\label{37+}
\begin{cases}
 &\dif \bar{x}(t)+(H^{-1}(t))^{2}\partial I_{\Xi}(\bar{x}(t))\dif t\ni \bar{f}(\bar{x}(t),  \bar{\mu}_{t})\dif t +\bar{g}(  \bar{x}(t),    \bar{\mu}_{t})\dif  B(t), t\in [t_{0}, T],\\
 & \bar{x}(t_{0})=H^{-1}(t)x_{0}.
 \end{cases}
 \end{equation}
where
\begin{align*}
&\bar{f}(\bar{x}(t),  \bar{\mu}_{t})=H^{-1}(t)(f(H(t)\bar{x}(t),  \bar{\mu}_{t})+H'(t)\bar{x}(t)),\\
&\bar{g}(\bar{x}(t),   \bar{\mu}_{t})=H^{-1}(t)(g(H(t)\bar{x}(t),  \bar{\mu}_{t})+H'(t)\bar{x}(t)).
\end{align*}
All assumptions in $(\mathrm{(A)}1)$ and $(\mathrm{A}2)$ are satisfied for coefficients of equation \eqref{37+}.  Consequently, using Theorem \ref{50}, Eq.\eqref{37+} admits a unique solution. Thus,  Eq.\eqref{37} has a solution.

 }
\end{exa}

\section{Stochastic principle of optimality} In this section, we will investigate optimal control for Eq.\eqref{40} below. The aim is  to show that the value function satisfies the dynamic  programming principle (DPP, for short).  Let $\mathbf{U}$ be a separable metric space. For a control process $u(\cdot,\cdot): [0,T]\times \Omega \rightarrow \mathbf{U},$ we consider the following stochastic controlled system:
\begin{align}\label{40}
\begin{cases}
 &\dif x(t)+ H(t)\partial \Pi(x(t))\dif t\ni f(x(t),  \mu_{t}, u(t))\dif t
   +g(  x(t),  \mu_{t}, u(t))\dif  B(t), t\in [s, T],\\
 & x(s)=x_{0},
 \end{cases}
 \end{align}
where $$f(\cdot, \cdot, \cdot ):  \mR^{m} \times \sP_{2}(\mR^{m})\times \mathbf{U} \rightarrow \mR^{m}, g(\cdot, \cdot, \cdot ):  \mR^{m} \times \sP_{2}(\mR^{m})\times \mathbf{U}\rightarrow \mR^{m\times d},$$
For the sake of simplicity, we assume $f(0, \delta_{0},u )=g(0,\delta_{0},u)=0$ and  make the following assumptions:
\begin{itemize}
\item[(A3)]The coefficients $f$ and  $g $ satisfy that  for some positive constant $L>0$ and all $x,y \in \mR^{d}, \mu, \nu \in \sP_{2}(\mR^{m}), u\in \mathbf{U}.$
\begin{align*}
|&f(x,\mu,u)-f(y,\nu,u)|+|g(x,\mu,u)-g(y,\nu,u)|
 \leq L(|x-y|+W_{2}(\mu, \nu)).
\end{align*}
From the above assumptions,  one can see that
$$|f(x,\mu,u)|\leq L(|x|+W_{2}(\mu, \delta_{0})), |g(x,\mu,u)|\leq L(|x|+W_{2}(\mu, \delta_{0})).$$
\end{itemize}

Define  the cost functional as follows:
\begin{align}\label{41}
J(s,x_{0};u)=\mE\bigg[\bigg.\int^{T}_{s}b(x^{s,x_{0},u}(t),  u(t))\dif t +\alpha(x^{s,x_{0},u}(T))\bigg]\bigg..
 \end{align}
where  $x^{s,x_{0},u}$ is the solution of Eq.\eqref{40} associated with $s, x_{0}, u$ and  $b: \mR^{m}\times \mathbf{U}\rightarrow \mR, \alpha: \mR^{m}\rightarrow \mR$ satisfy the following conditions:
\begin{itemize}
\item[$(\mathrm{A4})$]
For some positive constant $L>0$ and all $x,y \in \mR^{m}, \mu, \nu \in \sP_{2}(\mR^{m}), u\in \mathbf{U},$
\begin{align*}
|b(x,u)-b(y,u)|+|\alpha(x)-\alpha(y)|
 \leq L|x-y|,
\end{align*}
For the sake of simplicity, we assume $b(0, u )=\alpha(0)=0.$ This implies
$$ |b(x,u)|\leq L|x|, \alpha(x)\leq L|x|. $$
\end{itemize}

We  give the associated valued function as the infimum among all $u\in \sU[s,T]:$
\begin{align}\label{42}
V(s,x_{0}):=\inf_{u\in \sU[s,T]}J(s,x_{0};u), (s, x_{0})\in [0,T)\times \mR^{m},
 \end{align}
where $\sU[s, T]$ denotes the set of all the processes
$u(\cdot ,\cdot): \Omega \times [s, T]\rightarrow \mathbf{U}$  satisfying $$\mE\bigg[\bigg.\int^{T}_{s}|b(x^{s,x_{0},u}(r), u(r) )|\dif r+ |\alpha(x^{s,x_{0},u}(T) )|\bigg]\bigg.<\infty.$$

\bd \label{d1+}
If for any $(s,x_{0})\in [0,T)\times \mR^{m},$  it holds that
\begin{align}\label{d1++}
V(s,x_{0})=\inf_{u\in \sU[s,T]}\mE\bigg[\bigg.\int^{\tau}_{s}b(x^{s,x_{0},u}(t),  u(t))\dif t  +V(\tau, x^{s,x_{0},u}(\tau))\bigg]\bigg.,
 \end{align}
for any  $\tau \in [s,T].$  Then, we say that the value function $V$ satisfies the DDP.
\ed
In order to prove that the value function for Eq. (\ref{40}) fulfill DDP, we  consider the following penalized equation:
\begin{align}\label{43}
\begin{cases}
 & x^{\epsilon, s,x_{0},u}(t)+ \int^{t}_{s}H(r)\nabla \Pi_{\epsilon}(x^{\epsilon, s,x_{0},u}(r))\dif r
 =x_{0}+\int^{t}_{s}f(x^{\epsilon, s,x_{0},u}(r),  \mu_{t}, u(r))\dif r \\
   &\quad \quad\quad \quad\quad \quad\quad \quad\quad \quad\quad \quad\quad \quad\quad \quad+\int^{t}_{s}g(  x^{\epsilon, s,x_{0},u}(r),  \mu_{t}, u(r))\dif  B(r), t\in [s, T],\\
 & x^{\epsilon, s, x_{0}, u}(s)=x_{0},
 \end{cases}
 \end{align}
and  for any $ (t,x_{0})\in [0,T)\times \mR^{m}$   the penalized valued function is defined as follows:

\begin{align}\label{44}
V_{\epsilon}(s,x_{0})=\inf_{u\in \sU[s,T]}\mE\bigg[\bigg.\int^{T}_{s}b(x^{\epsilon, s,x_{0},u}(t), u(t))\dif t +\alpha(x^{\epsilon, s,x_{0},u}(T))\bigg]\bigg.,
 \end{align}
and set

\begin{align}\label{41b}
J_{\epsilon}(s,x_{0};u)=\mE\bigg[\bigg.\int^{T}_{s}b(x^{\epsilon, s,x_{0},u}(t),  u(t))\dif t +\alpha(x^{\epsilon, s,x_{0},u}(T))\bigg]\bigg..
 \end{align}

\bl \label{le4.1}
Assume $(\mathrm{A2})-(\mathrm{A4})$. Let $(x^{s, x_{0}, u}, k^{s, x_{0}, u})$ and $(x^{s', x_{0}', u}, k^{s', x_{0}', u})$ be the solutions of Eq.\eqref{40} corresponding to the initial date $(s, x_{0})$ and $(s', x_{0}')$ respectively. Then, we have the following estimates:
\begin{align}\label{45}
\mE[\sup_{s\leq r \leq T}|x^{s,x_{0},u}(r)|^{2}]
  \leq  C(a_{H},b_{H},x_{0}, L,T),
 \end{align}
and
\begin{align}\label{46}
\mE&[\sup_{s\vee s' \leq r \leq T}|x^{s,x_{0},u}(r)-x^{s',x_{0}',u}(r)|^{2}] \no\\
 & \leq C(a_{H},b_{H},x_{0}', L,T)\mE[|x_{0}-x_{0}'|^{2}] +  C(a_{H},b_{H},x_{0}', L,T)|s-s'|.
 \end{align}
\el
\begin{proof}
Assume that    $U^{s,x_{0},u},  U^{s',x_{0}',u}$ are two processes such that $$\dif k^{s, x_{0}, u}(t)= U^{s,x_{0},u}(t) \dif t, \dif k^{s', x_{0}', u}(t)= U^{s',x_{0}',u}(t) \dif t.$$
Using It\^{o}'s formula,  \eqref{45}   can be easily derived.   We only prove \eqref{46}.  Assume that $s\geq s'.$  For any $t\in [s, T],$ we have
\begin{align}\label{47}
  x^{s,x_{0},u}(t)+ \int^{t}_{s}H(r)U^{s,x_{0},u}(r)\dif r&= x^{s,x_{0},u}(s)+\int^{t}_{s}f(x^{ s,x_{0},u}(r),  \mu_{t}, u(r))\dif r\no \\
&\quad+\int^{t}_{s}g( x^{ s,x_{0},u}(r),  \mu_{t}, u(r))\dif  B(r), t\in [s, T],
 \end{align}
and
\begin{align}\label{48}
  x^{s',x'_{0},u}(t)+ \int^{t}_{s}H(r)U^{s',x'_{0},u}(r)\dif r&= x^{s,x'_{0},u}(s)+\int^{t}_{s}f(x^{ s',x_{0}',u}(r),  \mu'_{t}, u(r))\dif r\no \\
&\quad+\int^{t}_{s}g( x^{ s',x'_{0},u}(r),  \mu'_{t}, u(r))\dif  B(r), t\in [s, T].
 \end{align}
Set $N(s)=-\frac{\dif }{\dif s}(H^{-\frac{1}{2}}(s))=\frac{1}{2}H^{-\frac{3}{2}}(s)\frac{\dif H}{\dif s}(s).$
By $(\mathrm{A2}),$ we have
\begin{align}\label{49}
|N(s)|\leq \frac{1}{2}|H^{-\frac{3}{2}}(s)|\bigg|\bigg.\frac{\dif H}{\dif s}(s) \bigg|\bigg.\leq \frac{1}{2}a_{H}^{-\frac{3}{2}}M,
\end{align}
where $M:=\sup_{0\leq s\leq T}\big|\big.\frac{\dif H}{\dif s}(s) \big|\big..$
Set  $\hat{x}(t)=x^{s,x_{0},u}(t)-x^{s',x_{0}',u}(t),    \tilde{x}(t) =H^{-\frac{1}{2}}(t)\hat{x}(t).$
Then, $ \tilde{x}(\cdot)$  satisfies the following equation:
\begin{align}\label{50}
 \tilde{x}(t)&= \tilde{x} (s)+\int^{t}_{s}\hat{x}(r)\dif H^{-\frac{1}{2}}(r)+ \int^{t}_{s}H^{-\frac{1}{2}}(r)\dif\hat{x}(r)\no\\
& =\int^{t}_{t_{0}}R(r)\dif r + \int^{t}_{t_{0}}H^{-\frac{1}{2}}(r)(g( x^{ s,x_{0},u}(r),  \mu_{r}, u(r))\no\\
&\quad \quad \quad\quad \quad \quad\quad \quad \quad\quad \quad-g( x^{ s',x'_{0},u}(r),  \mu'_{r}, u(r)))\dif B(r).
\end{align}
where \begin{align*}
R(r)&:=\hat{x}(r)N(r)-H^{\frac{1}{2}}(r)[U^{s,x_{0},u}(r) - U^{s',x_{0}',u}(r)]\\
& \quad+H^{-\frac{1}{2}}(s)[f(x^{ s,x_{0},u}(r),  \mu_{r}, u(r))-f(x^{ s',x_{0}',u}(r),  \mu'_{r}, u(r))].
\end{align*}
Applying It\^{o}'s formula, we have
\begin{align}\label{51}
| \tilde{x}(t)|^{2}
& =| \tilde{x}(s)|^{2}+2\int^{t}_{s}\langle \tilde{x}(r), R(r)\rangle\dif r\no\\
&\quad + \int^{t}_{s} | H^{-\frac{1}{2}}(r)(g( x^{ s,x_{0},u}(r),  \mu_{r}, u(r))\no\\
&\quad \quad \quad\quad \quad \quad\quad \quad \quad\quad \quad-g( x^{ s',x'_{0},u}(r),  \mu'_{r}, u(r)))|^{2}\dif r\no\\
&\quad + 2\int^{t}_{s}\langle  \tilde{x}(r), H^{-\frac{1}{2}}(r)(g( x^{ s,x_{0},u}(r),  \mu_{r}, u(r))\no\\
&\quad \quad \quad\quad \quad \quad\quad \quad \quad\quad \quad-g( x^{ s',x'_{0},u}(r),  \mu'_{r}, u(r)))\rangle\dif B(r)\no\\
&\leq | \tilde{x}(s)|^{2}+(1+b_{H}^{\frac{1}{2}}a_{H}^{-\frac{3}{2}}M+4L^{2}a_{H}^{-1}b_{H}L^{2})\int^{t}_{s}\sup_{t_{0}\leq r\leq s}| \tilde{x}(r)|^{2}\dif s\no\\
&\quad + 4L^{2}a_{H}^{-1}b_{H}L^{2}\int^{t}_{s}W^{2}_{2}(\mu_{r}, \mu'_{r})\dif r\no\\
&\quad + 2\int^{t}_{s}\langle  \tilde{x}(r), H^{-\frac{1}{2}}(s)(g( x^{ s,x_{0},u}(r),  \mu_{r}, u(r))\no\\
&\quad \quad \quad\quad \quad \quad\quad \quad \quad\quad \quad-g( x^{ s',x'_{0},u}(r),  \mu'_{r}, u(r)))\rangle\dif B(r).
\end{align}
From BDG's inequality,  for any $l>0,$  we have
\begin{align*}
\mE&\sup_{s\leq r\leq t}\bigg|\bigg.\int^{r}_{s}\langle  \tilde{x}(u),  H^{-\frac{1}{2}}(u)(g( x^{ s,x_{0},u}(u),  \mu_{u}, u(u))\no\\
&\quad \quad \quad\quad \quad \quad\quad \quad \quad\quad \quad-g( x^{ s',x'_{0},u}(u),  \mu_{u}, u(u)))\rangle \dif  B(u)\bigg|\bigg.\\
&\leq 32\mE\bigg(\bigg.\int^{t}_{s}| \tilde{x}(s) |^{2}|H^{-\frac{1}{2}}(s)(g( x^{ s,x_{0},u}(r),  \mu'_{r}, u(r))\no\\
&\quad \quad \quad\quad \quad \quad\quad \quad \quad\quad \quad-g( x^{ s',x'_{0},u}(r),  \mu'_{r}, u(r)))|^{2} \dif r \bigg)\bigg.^{\frac{1}{2}}\\
& \leq 32\mE\bigg[\bigg. \sup_{s\leq r\leq t}| \tilde{x}(r)|\bigg(\bigg.\int^{t}_{s} |H^{-\frac{1}{2}}(s)(g( x^{ s,x_{0},u}(r),  \mu_{r}, u(r))\no\\
&\quad \quad \quad\quad \quad \quad\quad \quad \quad\quad \quad-g( x^{ s',x'_{0},u}(r),  \mu_{r}, u(r)))|^{2} \dif r \bigg)\bigg.^{\frac{1}{2}}\bigg]\bigg.\\
&\leq \frac{16}{l}\mE[\sup_{s\leq r\leq t}| \tilde{x}(r)|^{2}]  +16l  \mE\int^{t}_{s}|H^{-\frac{1}{2}}(r)(g( x^{ s,x_{0},u}(r),  \mu_{r}, u(r))\no\\
&\quad \quad \quad\quad \quad \quad\quad \quad \quad\quad \quad-g( x^{ s',x'_{0},u}(r),  \mu'_{r}, u(r)))|^{2}\dif r\\
&\leq  \frac{16}{l}\mE[\sup_{s\leq r\leq t}| \tilde{x}(r)|^{2}]  +l 64L^{2}a_{H}^{-1}b_{H}\mE\int^{t}_{s}| \tilde{x}(r)|^{2}\dif r,
\end{align*}
This together with \eqref{51} implies for  taking $l=2,$
\begin{align}\label{52}
&\mE[\sup_{s\leq r\leq t }| \tilde{x}(r)|^{2}]\no\\
&\leq  \mE[| \tilde{x}(s)|^{2}]
 +(2+2b_{H}^{\frac{1}{2}}a_{H}^{-\frac{3}{2}}M+(2^{11}+2^{5})L^{2}a_{H}^{-1}b_{H})\int^{t}_{s}\mE[\sup_{s\leq u\leq r }| \tilde{x}(u)|^{2}]\dif u.
\end{align}
Furthermore,
\begin{align}\label{53}
\mE[| \tilde{x}(s)|^{2}]
&\leq  3a_{H}^{-1}|x_{0}-x'_{0}|^{2}+  3a_{H}^{-1}\mE\bigg|\bigg.\int^{s}_{s'}f'(r)\dif r \bigg|\bigg.^{2}  +3a_{H}^{-1}\mE\int^{s}_{s'}|g'(r)|^{2}\dif  r \no\\
&\leq  3a_{H}^{-1}|x_{0}-x'_{0}|^{2}+  3a_{H}^{-1}(s-s')\mE\int^{s}_{s'}|f'(r)|^{2}\dif r
 +3a_{H}^{-1}\mE\int^{s}_{s'}|g'(r)|^{2}\dif r\no\\
&\leq 3a_{H}^{-1}|x_{0}-x'_{0}|^{2} + C(a_{H},b_{H},x'_{0}, L,T)(s-s'),
\end{align}
where $f'(r):= f(x'(r), \mu'_{r}), g'(r):= g(x'(r), \mu'_{r}).$
Combining \eqref{52} and \eqref{53},  we have
\begin{align}\label{54}
&\mE[\sup_{s\leq r\leq t }| \tilde{x}(r)|^{2}]
\leq  \mE[| \tilde{x}(s)|^{2}]\no\\
&\quad+(2+2b_{H}^{\frac{1}{2}}a_{H}^{-\frac{3}{2}}M+(2^{11}+2^{5})L^{2}a_{H}^{-1}b_{H})\int^{t}_{s}\mE[\sup_{s\leq u\leq r }| \tilde{x}(u)|^{2}]\dif u\no\\
&\leq (2b_{H}^{\frac{1}{2}}a_{H}^{-\frac{3}{2}}M+8L+\bar{\sigma}^{2}(6L^{2}+12L^{2}a_{H}^{-1}b_{H}))\int^{t}_{s}\mE[\sup_{s\leq u\leq r }| \tilde{x}(u)|^{2}]\dif u\no\\
&\quad + 3a_{H}^{-1}|x_{0}-x'_{0}|^{2} + C(a_{H},b_{H},x'_{0}, L,T, \overline{\sigma})(s-s').
\end{align}
Gronwall's inequality leads to
$$ \mE[\sup_{s\leq t\leq T }| \tilde{x}(t)|^{2}]\leq C(a_{H},b_{H},x'_{0}, L,T)|x_{0}-x'_{0}|^{2} +  C(a_{H},b_{H},x'_{0}, L,T)|s-s'|,   $$
The desired result is obtained.

\end{proof}

We can use techniques similar to those used in  Lemma \ref{le4.1} to give a priori estimate for solution $x^{\epsilon, s, x_{0},u}.$
\bl \label{le4.2}
Assume $(\mathrm{A2})-(\mathrm{A4})$. Let $x^{\epsilon, s, x_{0},u}$ be the solutions of Eq.\eqref{43} corresponding to the initial date $(s, x_{0}).$ Then, we have the following estimates:
\begin{align}\label{45l}
\mE[\sup_{s\leq r \leq T}|x^{\epsilon, s, x_{0},u}|^{2}]+\mE\int^{T}_{s}| \nabla\Pi_{\epsilon}(x^{\epsilon, s, x_{0},u}(s))|^{2}\dif s
  \leq  C(a_{H},b_{H},x_{0}, L,T).
 \end{align}

\el

The following lemma shows that $x^{\epsilon, s,x_{0},u}$ is a Cauchy sequence in $L^{2}(s, T; \mR^{m})$.

\bt \label{L4L}
Assume $(\mathrm{A}2)-(\mathrm{A}4).$  Then, we have the following estimates:
\begin{align}\label{+8+}
\mE[\sup_{t_{0}\leq t \leq T}|x^{\epsilon, s,x_{0},u}(t)-x^{\epsilon', s,x_{0},u}(t)|^{2}]\leq (a_{H},b_{H},x_{0}, L,T)(\epsilon+\epsilon').
\end{align}
\et
\begin{proof}
Set $N(s)=-\frac{\dif }{\dif s}(H^{-\frac{1}{2}}(s))=\frac{1}{2}H^{-\frac{3}{2}}(s)\frac{\dif H}{\dif s}(s).$
By $(\mathrm{A2}),$ we have
\begin{align}\label{9}
|N(s)|\leq \frac{1}{2}|H^{-\frac{3}{2}}(s)|\bigg|\bigg.\frac{\dif H}{\dif s}(s) \bigg|\bigg.\leq \frac{1}{2}a_{H}^{-\frac{3}{2}}M,
\end{align}
where $M:=\sup_{0\leq s\leq T}\big|\big.\frac{\dif H}{\dif s}(s) \big|\big..$
Set  $\hat{x}^{\epsilon, \epsilon'}(s)=x^{\epsilon, s,x_{0},u}(s)-x^{\epsilon', s,x_{0},u}(s),    \tilde{x}^{\epsilon, \epsilon'}(s) =H^{-\frac{1}{2}}(s)\hat{x}(s).$
Then, $\tilde{x}^{\epsilon, \epsilon'}(t)$  satisfies the following equation:
\begin{align}\label{10}
\tilde{x}^{\epsilon, \epsilon'}(t)&=-\int^{t}_{t_{0}}\hat{x}^{\epsilon, \epsilon'}(s)\dif H^{-\frac{1}{2}}(s)- \int^{t}_{t_{0}}H^{-\frac{1}{2}}(s)\dif\hat{x}^{\epsilon, \epsilon'}(s)\no\\
& =\int^{t}_{t_{0}}Q(s)\dif s + \int^{t}_{t_{0}}H^{-\frac{1}{2}}(s)(g( x^{\epsilon, s,x_{0},u}(s), \mu^{\epsilon}_{t})-g(  x^{\epsilon', s,x_{0},u}(s), \mu^{\epsilon'}_{s})\dif B(s)
\end{align}
where $Q(s):=\hat{x}^{\epsilon, \epsilon'}(s)N(s)-H^{\frac{1}{2}}(s)[\nabla \Pi_{\epsilon}(x^{\epsilon}(s))- \nabla \Pi_{\epsilon}(x^{\epsilon'}(s))]+H^{-\frac{1}{2}}(s)[f( x^{\epsilon, s,x_{0},u}(s), \mu^{\epsilon}_{s})-f( x^{\epsilon', s,x_{0},u}(s), \mu^{\epsilon'}_{s})].$
Applying It\^{o}'s formula and the property $(e)$ of $\nabla \Pi_{\epsilon}(\cdot),$ we have
\begin{align}\label{11}
|&\tilde{x}^{\epsilon, \epsilon'}(t)|^{2}
=2\int^{t}_{t_{0}}\langle\tilde{x}^{\epsilon, \epsilon'}(s), Q(s)\rangle\dif s\no\\
&\quad + \int^{t}_{t_{0}} | H^{-\frac{1}{2}}(s)(g( x^{\epsilon, s,x_{0},u}(s), \mu^{\epsilon}_{s})-g( x^{\epsilon', s,x_{0},u}(s), \mu^{\epsilon'}_{s})|^{2}\dif s\no\\
&\quad + 2\int^{t}_{t_{0}}\langle \tilde{x}^{\epsilon, \epsilon'}(s), H^{-\frac{1}{2}}(s)(g( x^{\epsilon, s,x_{0},u}(s), \mu^{\epsilon}_{t})-g( x^{\epsilon', s,x_{0},u}(s), \mu^{\epsilon'}_{s})\rangle\dif B(s)\no\\
&\leq C(a_{H},b_{H},\xi,t_{0}, L,T)\int^{t}_{t_{0}}\sup_{t_{0}\leq r\leq s }|\tilde{x}^{\epsilon, \epsilon'}(r)|^{2}\dif s\no\\
&\quad +(\epsilon+\epsilon')\int^{t}_{t_{0}}|\nabla \Pi_{\epsilon}(x^{\epsilon}(s))|^{2}\dif s+ (\epsilon+\epsilon')\int^{t}_{t_{0}}|\nabla \Pi_{\epsilon'}(x^{\epsilon'}(s))|^{2}\dif s\no\\
&\quad + 2\int^{t}_{t_{0}}\langle \tilde{x}^{\epsilon, \epsilon'}(s), H^{-\frac{1}{2}}(s)(g( x^{\epsilon, s,x_{0},u}(s), \mu^{\epsilon}_{s})-g( x^{\epsilon', s,x_{0},u}(s), \mu^{\epsilon'}_{s})\rangle\dif B(s).
\end{align}
From BDG's inequality,  for any $l>0,$  we have
\begin{align*}
\mE&\sup_{0\leq r\leq t}\bigg|\bigg.\int^{r}_{t_{0}}\langle \tilde{x}^{\epsilon, \epsilon'}(s),  H^{-\frac{1}{2}}(s)(g( x^{\epsilon, s,x_{0},u}(s), \mu^{\epsilon}_{s})-g( x^{\epsilon', s,x_{0},u}(s), \mu^{\epsilon'}_{s})\rangle \dif  B(s)\bigg|\bigg.\\
& \leq \mE\bigg(\bigg.\int^{t}_{t_{0}}|\tilde{x}^{\epsilon, \epsilon'}(s) |H^{-\frac{1}{2}}(s)(g( x^{\epsilon, s,x_{0},u}(s), \mu^{\epsilon}_{s})-g( x^{\epsilon', s,x_{0},u}(s), \mu^{\epsilon'}_{s})|^{2} \dif s \bigg)\bigg.^{\frac{1}{2}}\\
& \leq \mE\bigg[\bigg. \sup_{t_0\leq s\leq t}|\tilde{x}^{\epsilon, \epsilon'}(s)|\bigg(\bigg.\int^{t}_{t_{0}} |H^{-\frac{1}{2}}(s)(g( x^{\epsilon, s,x_{0},u}(s), \mu^{\epsilon}_{s})-g( x^{\epsilon', s,x_{0},u}(s), \mu^{\epsilon'}_{s}|^{2} \dif s \bigg)\bigg.^{\frac{1}{2}}\bigg]\bigg.\\
& \leq \frac{1}{l}\mE[\sup_{t_0\leq s\leq t}|\tilde{x}^{\epsilon, \epsilon'}(s)|^{2}]  +l  \mE\int^{t}_{t_{0}}|H^{-\frac{1}{2}}(s)(g( x^{\epsilon, s,x_{0},u}(s), \mu^{\epsilon}_{s})-g( x^{\epsilon', s,x_{0},u}(s), \mu^{\epsilon'}_{s})|^{2}\dif s\\
&\leq  \frac{1}{l}\mE[\sup_{t_0\leq s\leq t}|\tilde{x}^{\epsilon, \epsilon'}(s)|^{2}]  +l C(a_{H},b_{H},\xi,t_{0}, L,T) \mE\int^{t}_{t_{0}}| \tilde{x}^{\epsilon, \epsilon'}(s)|^{2}\dif s,
\end{align*}
This together with \eqref{11} implies for  taking $l=2,$
\begin{align*}
&\mE[\sup_{t_{0}\leq t\leq T }|\tilde{x}^{\epsilon, \epsilon'}(s)|^{2}]\\
&\leq C(a_{H},b_{H},\xi,t_{0}, L,T) \int^{t}_{t_{0}}\mE[\sup_{t_{0}\leq r\leq s }|\tilde{x}^{\epsilon, \epsilon'}(r)|^{2}]\dif s\no\\
&\quad +2(\epsilon+\epsilon')\mE\int^{t}_{t_{0}}|\nabla \Pi_{\epsilon}(x^{\epsilon}(s))|^{2}\dif s+ 2(\epsilon+\epsilon')\mE\int^{t}_{t_{0}}|\nabla \Pi_{\epsilon'}(x^{\epsilon'}(s))|^{2}\dif s\no\\
&\leq C(a_{H},b_{H},\xi,t_{0}, L,T)\int^{t}_{t_{0}}\mE[\sup_{t_{0}\leq r\leq s }|\tilde{x}^{\epsilon, \epsilon'}(r)|^{2}]\dif s\\
&\quad  +(\epsilon+\epsilon')C(a_{H},b_{H},x_{0}, L,T).
\end{align*}
Gronwall's inequality leads to
$$ \mE[\sup_{t_{0}\leq t\leq T }|\tilde{x}^{\epsilon, \epsilon'}(s)|^{2}]\leq C(a_{H},b_{H},x_{0}, L,T)(\epsilon+\epsilon').   $$
The proof is therefore complete.

\end{proof}

\bl\label{l4.4}
Let $(\mathrm{A}2)-(\mathrm{A}4)$ hold.   The following equality holds for a subsequence of $\epsilon,$ which we still denote by $\epsilon, $
$$\lim_{\epsilon \rightarrow 0}\mE\int^{T}_{s}|x^{\epsilon, s,x_{0},u}(r)- x^{ s,x_{0},u}(r)|^{2}\dif r=0. $$

\el
\begin{proof}
By Lemma \ref{L4L}, we know that there exists a stochastic process $\check{x}^{ s,x_{0},u}\in L^{2}(s, t; \mR^{m})$ such that $\lim_{\epsilon\rightarrow 0}x^{\epsilon, s,x_{0},u}=\check{x}^{ s,x_{0},u} $ in $L^{2}(s, t; \mR^{m}).$   Moreover, by   Lemma \ref{le4.2}, there exists a stochastic process $\check{U}^{ s,x_{0},u}\in L^{2}(s, T; \mR^{m})$ such that $\nabla\Pi_{\epsilon}(x^{\epsilon, s,x_{0},u}(\cdot))\Rightarrow \check{U}^{ s,x_{0},u}(\cdot) $ in $L^{2}(s, T; \mR^{m}), \epsilon\rightarrow 0. $   As a consequence,
$$ \lim_{\epsilon\rightarrow 0}\int^{T}_{t_{0}}H(s)\nabla\Pi_{\epsilon}(x^{\epsilon, s,x_{0},u}(s))\dif s =\lim_{\epsilon\rightarrow 0}\int^{T}_{t_{0}}H(s)\check{U}^{ s,x_{0},u}(s)\dif s. $$
Therefore, we can pass to the limit in the approximating equation \eqref{43} and  obtain that
\begin{align*}
\begin{cases}
 &\check{x}^{ s,x_{0},u}(t)+ \int^{t}_{s}H(r)\check{U}^{ s,x_{0},u}(r)\dif r=  \int^{t}_{s}f(\check{x}^{ s,x_{0},u}(r),  \check{\mu}^{ s,x_{0},u}_{r}, u(r))\dif r\\
   &\quad \quad \quad \quad \quad \quad\quad \quad \quad\quad \quad \quad+\int^{t}_{s}g(  \check{x}^{ s,x_{0},u}(r),  \check{\mu}^{ s,x_{0},u}_{r}, u(r))\dif  B(r), t\in [s, T],\\
 & \check{x}(s)=x_{0},
 \end{cases}
\end{align*}
where $\check{\mu}^{ s,x_{0},u}_{t}$ is the distribution of $\check{x}^{ s,x_{0},u}(t).$
By $(b),$  we have $$\nabla \Pi_{\epsilon}(x^{\epsilon, s,x_{0},u}(t))\in \partial \Pi(J_{\epsilon}(x^{\epsilon, s,x_{0},u}(t))).$$  Thus,

\begin{align}\label{90}
\mE&\int^{T}_{s}\langle\nabla \Pi_{\epsilon}(x^{\epsilon, s,x_{0},u}(r)), \nu(r) -J_{\epsilon}(x^{\epsilon, s,x_{0},u}(t)))  \rangle \dif r +  \mE\int^{T}_{s} \Pi(J_{\epsilon}(x^{\epsilon, s,x_{0},u}(r))) \dif r\no\\
&\leq   \mE\int^{T}_{s} \Pi(\nu(r)) \dif r
\end{align}
for any square integrable stochastic process $\nu.$
From the definition of $J_{\epsilon}(\cdot)$  and Lemma \ref{le4.2}, we know that

\begin{align*}
 &\lim_{\epsilon\rightarrow 0}\int^{T}_{s}|J_{\epsilon}(x^{\epsilon, s,x_{0},u}(t))- \check{x}^{ s,x_{0},u}(r)|^{2}\dif r\\
 &=2\epsilon^{2}\limsup_{\epsilon \rightarrow 0}\mE\int^{T}_{s}|\nabla\Pi_{\epsilon}(x^{\epsilon, s,x_{0},u}(r))|^{2}\dif r+2\limsup_{\epsilon \rightarrow 0}\mE\int^{T}_{s}|x^{\epsilon, s,x_{0},u}(r)-\check{x}^{ s,x_{0},u}(r)|^{2}\dif r\\
 &=0.
 \end{align*}
By taking $\liminf_{\epsilon\rightarrow 0}$ on both side of the above inequality \eqref{90}, we have $U(s)\in \partial \Pi(x(s)).$
Thus $\check{x}^{ s,x_{0},u}$ be the solutions of Eq.\eqref{40} corresponding to the initial date $(s, x_{0}).$  By the uniqueness, we know that $\check{x}^{ s,x_{0},u}(t)= x^{ s,x_{0},u}(t), a.e.$  Then,

\end{proof}

\bl\label{h+}
Assume $(\mathrm{A2})-(\mathrm{A4}).$   The following estimates hold:
\begin{align}
&|V(s, x_{0})| \leq  C(a_{H},b_{H},x_{0}, L,T),\label{b+}\\
&|V(s, x_{0})- V(s', x'_{0})| \nonumber\\
&\leq C(a_{H},b_{H},x'_{0}, L,T)|x_{0}-x'_{0}|+ C(a_{H},b_{H},x'_{0}, L,T)|s-s'|^{\frac{1}{2}}.\label{b++}
\end{align}

\el
\begin{proof}
By $(\mathrm{A3}),$ we have
\begin{align*}
|J(t, x_{0})|&=\bigg|\bigg. \mE\bigg[\bigg. \int^{T}_{s}b(t, x^{s, x_{0}, u}(t), u(t))\dif t + \alpha(x^{s, x_{0}, u}(T))  \bigg]\bigg.    \bigg|\bigg.\\
& \leq C(T,L)\mE(1+\sup_{s\leq t\leq T}|x^{s, x_{0}, u}(t)|)\\
&\leq C(a_{H},b_{H},x_{0}, L,T).
\end{align*}
Furthermore, for the sake of simplicity, assume $s\geq s'.$   Then,  we have
\begin{align}\label{c++}
|&J(s, x_{0})- J(s', x'_{0})|=\bigg|\bigg. \mE\bigg[\bigg. \int^{T}_{s}b(t, x^{s, x_{0}, u}(t),  u(t))\dif t + \alpha(x^{s, x_{0}, u}(T))  \bigg]\bigg.\no\\
&\quad - \mE\bigg[\bigg. \int^{T}_{s'}b(t, x^{s',  x'_{0}, u}(t),  u(t))\dif t + \alpha( x^{s', x'_{0}, u}(T))  \bigg]\bigg.   \bigg|\bigg.\no\\
& \leq \mE\bigg|\bigg.\int^{T}_{s}b(t, x^{s, x_{0}, u}(t),  u(t))-\int^{T}_{s'}b(t, x^{s', x'_{0}, u}(t), u(t))\dif t\bigg|\bigg. \no\\
&\quad + \mE[|\alpha( x^{s', x'_{0}, u}(T))- \alpha( x^{s', x'_{0}, u}(T))   |]\no\\
&\quad + \mE[|\alpha( x^{s, x_{0}, u}(T))- \alpha( x^{s', x'_{0}, u}(T))   |]\no\\
& \leq C( L, T)\mE[\sup_{s\leq t\leq T}|x^{s, x_{0}, u}(t)- x^{s', x'_{0}, u}(t)|]\no\\
&\quad + C( L, T)(s-s')+ C( L, T)(s-s')\mE[\sup_{s\leq t\leq T}| x^{s', x'_{0}, u}(t)|]\no\\
& \leq C(a_{H},b_{H},x'_{0}, L,T)|x_{0}-x'_{0}|+ C(a_{H},b_{H},x'_{0}, L,T)|s-s'|^{\frac{1}{2}}.
\end{align}
\end{proof}

The following lemma about the DDP of Eq.\eqref{43} is a straightforward result  of Theorem 4.2 in \cite{LC}.
\bl
 Assume $(\mathrm{A2})-(\mathrm{A4})$.   Let $x^{\epsilon, s,x_{0},u}$ be the solution of \eqref{43}. Then, for  any $ (t,x_{0})\in [0,T)\times \mR^{m},$   it holds that
 \begin{align}\label{55}
V_{\epsilon}(s,x_{0})=\inf_{u\in \sU[s,T]}\mE\bigg[\bigg.\int^{\tau}_{s}b(x^{\epsilon, s,x_{0},u}(t),  u(t))\dif t  +V(\tau, x^{\epsilon, s,x_{0},u}(\tau))\bigg]\bigg.,
 \end{align}
 for every stopping time $\tau\in  [s, T].$
 \el

Now,  we are going to prove that Eq.\eqref{40} satisfies the DDP
\bp
Assume $(\mathrm{A2})-(\mathrm{A4})$. Then it holds that:
\begin{align}\label{56}
&|V_{\epsilon}(s,x_{0})-V(s,x_{0})|\leq  C(a_{H},b_{H},x_{0}, L,T)\epsilon^{\frac{1}{2}}.
 \end{align}
 \ep
\begin{proof}
By Lemma \ref{l4.4}, taking the limit in \eqref{+8+}, we deduce that
$$\mE[\sup_{s\leq r\leq T}|x^{\epsilon, s,x_{0},u}(r)- x^{ s,x_{0},u}(r)|^{2}]\leq C(a_{H},b_{H},x_{0}, L,T)\epsilon.$$
With the same argument as in the proof of  Lemma \ref{h+}, we obtain
\begin{align*}
|J_{\epsilon}(s,x_{0},u)-J(s,x_{0},u)|\leq C(a_{H},b_{H},x_{0}, L,T)\epsilon^{\frac{1}{2}}.
 \end{align*}
Consequently, we obtain
$$|V_{\epsilon}(s,x_{0})-V(s,x_{0})|\leq  \sup_{u\in \sU[s, T]}|J_{\epsilon}(s,x_{0},u)-J(s,x_{0},u)| \leq  C(a_{H},b_{H},x_{0}, L,T)\epsilon^{\frac{1}{2}}.      $$
The second   conclusion follows.

\end{proof}

From the above lemma, the following result immediately holds.

\bl
$V_{\epsilon}(\cdot, \cdot)$ is uniformly convergent on compacts  to the value function $V(\cdot, \cdot)$ on $[0,T)\times \overline{D(\Pi)}.$
\el

Now, we state our main result.

\bt
Under the assumptions $(\mathrm{A2})-(\mathrm{A4}),$   the value function $V$ satisfies the DPP \eqref{d1++}.

\et
\begin{proof}
Taking any $(t,x_{0})\in [0,T)\times \overline{D(\Pi)},$  we have, for any $\epsilon >0, u\in \sU[s, T], \tau \in [s, T],$ and $R>0,$
\begin{align}\label{33+}
&\mE[|V_{\epsilon}(\tau,x^{\epsilon, s,x_{0},u}(\tau))-V(\tau,x^{ s,x_{0},u}(\tau))|]\no\\
&\leq\mE[|V_{\epsilon}(\tau,x^{\epsilon, s,x_{0},u}(\tau))-V_{\epsilon}(\tau,x^{ s,x_{0},u}(\tau))|]\no\\
&\quad +\mE[|V_{\epsilon}(\tau,x^{ s,x_{0},u}(\tau))-V(\tau,x^{ s,x_{0},u}(\tau))|(1_{A_{1}}+1_{A_{2}})]\no\\
& \leq \mE[|V_{\epsilon}(\tau,x^{\epsilon, s,x_{0},u}(\tau))-V_{\epsilon}(\tau,x^{ s,x_{0},u}(\tau))|] +\sup_{(t,y)\in B}\mE[|V_{\epsilon}(t,y)-V(t,y)|]\no\\
&\quad +\mE[|V_{\epsilon}(\tau,x^{ s,x_{0},u}(\tau))-V(\tau,x^{ s,x_{0},u}(\tau))|1_{A_{2}}],
 \end{align}
where
$$A_{1}:= \{\omega: |x^{ s,x_{0},u}(\tau)|\leq R  \}, A_{2}:= \{\omega: |x^{ s,x_{0},u}(\tau)|> R  \},  $$
and
$$B:= [0, T]\times (\overline{D(\Pi)}\cap \overline{B(0, R)}).$$
With the same argument as in the proof of  Lemma \ref{h+},  we have
\begin{align}\label{57}
|J_{\epsilon}(\tau, x^{ \epsilon, s,x_{0},u}(\tau))- J_{\epsilon}(\tau, x^{ s,x_{0},u}(\tau))|\leq C(a_{H},b_{H},\xi, L,T)\epsilon^{\frac{1}{2}}.
\end{align}
Now, we look at the term $\mE[|V_{\epsilon}(\tau,x^{ s,x_{0},u}(\tau))-V(\tau,x^{ s,x_{0},u}(\tau))|1_{A_{2}}].$  By definition of $A_{2},$ we have
\begin{align}\label{58}
\mE&[|V_{\epsilon}(\tau,x^{ s,x_{0},u}(\tau))-V(\tau,x^{ s,x_{0},u}(\tau))|1_{A_{2}}]\leq (\mE[|V_{\epsilon}(\tau,x^{ s,x_{0},u}(\tau))-V(\tau,x^{ s,x_{0},u}(\tau))|^{2}])^{\frac{1}{2}}(\mE1_{A_{2}})^{\frac{1}{2}}\no\\
&\leq \sqrt{2}(\mE[|V_{\epsilon}(\tau,x^{ s,x_{0},u}(\tau))|^{2}]+\mE[|V(\tau,x^{ s,x_{0},u}(\tau))|^{2}])^{\frac{1}{2}}|]\frac{(\mE| x^{ s,x_{0},u}(\tau)|^{2})^{\frac{1}{2}}}{R}\no\\
&\leq \frac{C(a_{H},b_{H},x_{0}, L,T)}{R}.
\end{align}
Combining \eqref{33+}, \eqref{57} and \eqref{58},  we have
\begin{align}\label{59}
&\mE[|V_{\epsilon}(\tau,x^{\epsilon, s,x_{0},u}(\tau))-V(\tau,x^{ s,x_{0},u}(\tau))|]\no\\
&\leq \mE[|V_{\epsilon}(\tau,x^{\epsilon, s,x_{0},u}(\tau))-V_{\epsilon}(\tau,x^{ s,x_{0},u}(\tau))|] +\sup_{(t,y)\in B}\mE[|V_{\epsilon}(t,y)-V(t,y)|]\no\\
&\quad+\frac{C(a_{H},b_{H},x_{0}, L,T)}{R}\no\\
&\leq  C(a_{H},b_{H},x_{0}, L,T)\epsilon^{\frac{1}{2}}+\sup_{(t,y)\in B}\mE[|V_{\epsilon}(t,y)-V(t,y)|]\no\\
&\quad+\frac{C(a_{H},b_{H},x_{0}, L,T)}{R}.
 \end{align}
In addition,
\begin{align*}
V(t, \xi)&\leq V_{\epsilon}(t, x_{0})+|V_{\epsilon}(t, x_{0})-V(t, x_{0})|\no\\
& \leq \mE\bigg[\bigg.\int^{\tau}_{s}b(x^{\epsilon, s,x_{0},u}(r), x^{\epsilon, s,x_{0},u}_{r}, u(r))\dif r  +V(\tau, x^{\epsilon, s,x_{0},u}(\tau))\bigg]\bigg.+|V_{\epsilon}(t, x_{0})-V(t, x_{0})|\no\\
& \leq \mE\bigg[\bigg.\int^{\tau}_{s}b(x^{ s,x_{0},u}(r), u(r))\dif r  +V(\tau, x^{ s,x_{0},u}(\tau))\bigg]\bigg.+|V_{\epsilon}(t, x_{0})-V(t,x_{0})|\no\\
&\quad +\mE\bigg[\bigg.\int^{\tau}_{s}b(x^{\epsilon, s,x_{0},u}(r),  u(r))\dif r  -\int^{\tau}_{s}b(x^{ s,x_{0},u}(r),  u(r))\dif r\bigg]\bigg.\no\\
& \quad +|V_{\epsilon}(\tau, x^{\epsilon, s,x_{0},u}(\tau))-V(\tau, x^{ s,x_{0},u}(\tau))|\no\\
& \leq \mE\bigg[\bigg.\int^{\tau}_{s}b(x^{ s,x_{0},u}(r),  u(r))\dif r  +V(\tau, x^{ s,x_{0},u}(\tau))\bigg]\bigg.\no\\
&\quad +  C(a_{H},b_{H},x_{0}, L,T)\epsilon^{\frac{1}{2}}+\sup_{(t,y)\in B}\mE[|V_{\epsilon}(t,y)-V(t,y)|]\no\\
&\quad+\frac{C(a_{H},b_{H},\xi, L,T)}{R}.
\end{align*}
Passing to the limit for $\epsilon\rightarrow 0$ and $R\rightarrow \infty,$ we obtain
\begin{align*}
V(t, x_{0})
\leq  \mE\bigg[\bigg.\int^{\tau}_{s}b(x^{ s,x_{0},u}(r),  u(r))\dif r  +V(\tau, x^{ s,x_{0},u}(\tau))\bigg]\bigg..
\end{align*}
Conversely, for any $\delta> 0,$ since $V_{\epsilon}$ satisfies the DPP, there exists $u_{\delta}\in \sU[t,T]$ such that
$$V_{\epsilon}(t,x)+\frac{\delta}{2}\geq \mE\bigg[\bigg.\int^{\tau}_{s}b(x^{\epsilon, s,\xi,u_{\delta}}(r),  u_{\delta}(r))\dif r  +V_{\epsilon}(\tau, x^{\epsilon, s,x_{0},u_{\delta}}(\tau))\bigg]\bigg.. $$
Using the above inequality,  we have
\begin{align*}
&V(t,x)+\delta\geq V_{\epsilon}(t,x)+\delta-|V_{\epsilon}(t,x)-V(t,x)|\\
&\geq \mE\bigg[\bigg.\int^{\tau}_{s}b(x^{\epsilon, s,x_{0},u_{\delta}}(r),  u_{\delta}(r))\dif r  +V_{\epsilon}(\tau, x^{\epsilon, s,x_{0},u_{\delta}}(\tau))\bigg]\bigg.+\frac{\delta}{2}-|V_{\epsilon}(t,x)-V(t,x)|\no\\
&\geq \mE\bigg[\bigg.\int^{\tau}_{s}b(x^{ s,x_{0},u_{\delta}}(r), u_{\delta}(r))\dif r  +V(\tau, x^{ s,x_{0},u_{\delta}}(\tau))\bigg]\bigg.+\frac{\delta}{2}-|V_{\epsilon}(t, x_{0})-V(t, x_{0})|\no\\
&\quad -\mE\bigg[\bigg.\int^{\tau}_{s}b(x^{\epsilon, s,x_{0},u_{\delta}}(r),  u_{\delta}(r))\dif r  -\int^{\tau}_{s}b(x^{ s,x_{0},u_{\delta}}(r), u_{\delta}(r))\dif r\bigg]\bigg.\no\\
&\quad -|V_{\epsilon}(\tau, x^{\epsilon, s,x_{0},u_{\delta}}(\tau))-V(\tau, x^{ s,x_{0},u_{\delta}}(\tau))|\\
&\geq \mE\bigg[\bigg.\int^{\tau}_{s}b(x^{ s,x_{0},u_{\delta}}(r), x^{ s,x_{0},u_{\delta}}_{r}, u_{\delta}(r))\dif r  +V(\tau, x^{ s,x_{0},u_{\delta}}(\tau))\bigg]\bigg.+\frac{\delta}{2}-|V_{\epsilon}(t, x_{0})-V(t, x_{0})|\no\\
&\quad -  C(a_{H},b_{H},x_{0}, L,T)\epsilon^{\frac{1}{2}}-\sup_{(t,y)\in B}\mE[|V_{\epsilon}(t,y)-V(t,y)|]-\frac{C(a_{H},b_{H},x_{0}, L,T)}{R}\no\\
&\geq \mE\bigg[\bigg.\int^{\tau}_{s}b(x^{ s,\xi,u_{\delta}}(r), x^{ s,x_{0},u_{\delta}}_{r}, u_{\delta}(r))\dif r  +V(\tau, x^{ s,x_{0},u_{\delta}}(\tau))\bigg]\bigg.,
\end{align*}
where the last inequality follows by letting $\epsilon\rightarrow 0$ and $R\rightarrow \infty.$ The proof is complete.

\end{proof}

\section*{Funding}
This research was supported by the National Natural Science Foundation of China (Grant nos. 61876192, 11626236), the Fundamental Research Funds for the Central Universities of South-Central University for Nationalities (Grant nos. CZY15017, KTZ20051, CZT20020).
\section*{Availability of data and materials}
\begin{center}
Not applicable.
\end{center}

\section*{Competing interests}
\begin{center}
The author declare they have no competing interests.
\end{center}

\section*{Authors' contributions}
All authors conceived of the study and participated in its design and coordination. All authors read and approved the final manuscript.

\end{document}